%% file: slcgenproj.tex
\documentclass{amsart}
\usepackage{amssymb, latexsym, enumerate, amscd, amsthm, amsmath}

\input{hart}

\theoremstyle{plain}
\newtheorem{thm}{Theorem}[section]
\newtheorem{cor}[thm]{Corollary}
\newtheorem{prop}[thm]{Proposition}
\newtheorem{lemma}[thm]{Lemma}
\newtheorem*{mthm}{Main Theorem}
\newtheorem*{thm42}{Theorem 4.2}
\newtheorem*{cor47}{Corollary 4.7}

\theoremstyle{definition}
\newtheorem{defn}[thm]{Definition}
\newtheorem{ex}[thm]{Example}

\theoremstyle{remark}

\begin{document}

\title[Generic Projection Hypersurfaces]{Singularities of Generic Projection Hypersurfaces}

\author[D. C. Doherty]{Davis C. Doherty}

\begin{abstract}
Linearly projecting smooth projective varieties provides a method of obtaining hypersurfaces birational to the original varieties.  We show that in low dimension, the resulting hypersurfaces only have Du Bois singularities.  Moreover, we conclude that these Du Bois singularities are in fact semi log canonical.  However, we demonstrate the existence of counterexamples in high dimension -- the generic linear projection of certain varieties of dimension 30 or higher is neither semi log canonical nor Du Bois.
\end{abstract}

\maketitle

\section{Introduction}

The study of smooth projective curves over $\C$ is greatly simplified by the classical result that every such curve is birational to a plane curve having only nodal singularities, with the birational morphism given by a generic linear projection.  Applying an analogous technique to higher-dimensional varieties seems a natural extension -- can we draw any conclusions in this more general setting?  Joel Roberts (\cite{Rob71}, \cite{Rob75}) provides a useful starting point for understanding the singularities introduced by linearly projecting varieties.  Applying the more recent machinery of birational geometry to his work offers the possibility of finding a useful answer to this question.

Our approach is as follows: start with a smooth projective variety, obtain a birational projective hypersurface by taking a (generic) linear projection, and classify the singularities of the resulting hypersurface.  The principal challenge is in determining the type of singularities introduced by the projection -- in order for this approach to be useful, we would like to constrain these to a reasonably nice class.  To that end, we initially propose that the appropriate class of singularities is \textit{semi log canonical} (slc) -- these are the singularities appearing on the boundaries of moduli spaces for most higher-dimensional moduli problems. 

Unfortunately, determining semi log canonicity is highly nontrivial even for hypersurfaces.  To that end, the somewhat obscure class of Du Bois singularities serves as a useful tool, as we establish in the following.

\begin{thm42} Let $X$ be an $S_2$ scheme which is semismooth in codimension one, and assume that $K_{X}$ is Cartier and $X$ has Du Bois singularities.  Then $X$ is semi log canonical.
\end{thm42}

Applying a variety of techniques to the results of \cite{Rob75}, combined with the above, leads to the following extension of the classical result on curves.

\begin{mthm} Let $Y\subset\P^N$ be a smooth projective variety of dimension $n$, $n\leq5$, embedded via the $d$-uple embedding with $d\geq 3n$.  Let $X\subset\P^{n+1}$ be the image of $Y$ under a generic projection $\pi:Y\to\P^{n+1}$.  Then $X$ has semi log canonical singularities.
\end{mthm}

Though the dimensional restrictions in the Main Theorem seem arbitrary, the statement cannot be generalized under these assumptions: counterexamples exist in high dimensions.

\begin{cor47} Let $X\subset\P^{31}$ be a generic projection hypersurface obtained via $\pi:Y\to\P^{31}$, where $\OM^1_Y$ is nef.  Then $X$ is not semi log canonical.
\end{cor47}

An interesting fact to note is that the only known counterexamples all have $\OM^1_Y$ nef -- in particular, the originating smooth variety $Y$ is a minimal model.  This leaves open the possibility for a qualified version of the original hypothesis regarding the introduced singularities to hold.

I would like to express my gratitude to Karl Schwede, whose help with Du Bois singularities was instrumental, and Rob Lazarsfeld, who pointed out that the counterexamples of Corollary 4.7 should exist.

\section{Definitions and Conventions}

Throughout this work we assume all schemes to be separated of finite type over $\C$.  A variety is a reduced and separated scheme of finite type over $\C$.  A scheme is \textit{Gorenstein} if all its local rings are Gorenstein rings.  Similarly, scheme is $S_2$ if all its local rings satisfy Serre's $S_2$ condition.  A scheme (or more generally, an algebraic space) $X$ of dimension $n$ is \textit{semismooth} if every closed point is either smooth, a double normal crossing point -- analytically isomorphic to $\C[x_0,\ldots,x_n]/(x_0 x_1)$ -- or a pinch point -- analytically isomorphic to $\C[x_0,\ldots,x_n]/(x_0 x_1^2-x_2^2)$ (note that a semismooth scheme is Gorenstein).  Its \textit{double locus} $D_X$ is the codimension one subscheme of non-smooth points.  A proper birational map $f:Y\to X$ is called a \emph{semiresolution} if $Y$ is semismooth, no component of its double locus $D_Y$ is in the exceptional locus of $f$, and there is a codimension two subset $S\subset X$ such that the restriction map $f^{-1}(X\setminus S)\to X\setminus S$ is an isomorphism.  $f$ is a \emph{good semiresolution} if, in addition, $E\cup D_Y$ is a simple normal crossing divisor (where $E$ is the exceptional divisor of $f$).  By \cite[Prop 4.2]{K90} semiresolutions always exist, as long as we're willing to work in the category of algebraic spaces.

Now suppose $X$ is a reduced $S_2$ scheme which is semismooth in codimension one.  We say that $X$ has \emph{semi log canonical} singularities if $K_X$ is $\Q$-Cartier and there is a good semiresolution $f:Y\to X$ such that
\[
K_Y\equiv f^{*}K_X+\sum_{i}a_i E_i,
\]
with all $a_i\geq -1$, where $E_i$ are the exceptional divisors.

\subsection{Du Bois Singularities}

Denote by $D_{\text{filt}}(X)$ the bounded, filtered derived category of sheaves on a scheme $X$ with coherent cohomology; this is the only derived category we shall consider, so this abbreviated notation will be unambiguous.  We denote quasi-isomorphisms by $\qis$.

Philippe Du Bois (\cite{DB81}) demonstrated the existence of an object $\DB{\cdot}{X}\in D_{\text{filt}}(X)$ (the \emph{Du Bois complex}) with the following properties:
\begin{enumerate}[\upshape (i)]
\item If $\OM_{X}^{\cdot}$ denotes the usual De Rham complex with the ``filtration b\^ete'', then there is a natural morphism $\OM_{X}^{\cdot}\to\DB{\cdot}{X}$ in the filtered derived category (i.e., respecting the filtrations); if $X$ is smooth, this map is a quasi-isomorphism.

\item The Du Bois complex is local in the \'etale topology -- if $U$ is an \'etale open set, then $\DB{\cdot}{X}|_{U}\qis\DB{\cdot}{U}$.

\item If $f:Y\to X$ is a proper morphism, then there is a natural morphism $f^{*}:\DB{\cdot}{X}\to Rf_{*}\DB{\cdot}{Y}$ in the filtered derived category.

\item Let $f:Y\to X$ be proper, and assume that $f$ is an isomorphism outside a closed subscheme $\Sigma\subset X$, with reduced preimage $f^{-1}(\Sigma)=E$.  Then there exists an exact triangle
\[
\DB{\cdot}{X}\to\DB{\cdot}{\Sigma}\oplus Rf_{*}\DB{\cdot}{Y}\to Rf_{*}\DB{\cdot}{E}\xrightarrow{\text{\upshape +1}}.
\]
\end{enumerate}
Set $\DB{0}{X}=\Gr_{\text{filt}}^{0}\DB{\cdot}{X}$.  There is a natural morphism $\sO_X\to\DB{0}{X}$; we say $X$ has \emph{Du Bois singularities} if this is a quasi-isomorphism.

\section{Identifying Du Bois Singularities}

Determining when a scheme has Du Bois singularities is a difficult problem.  In this section, we review some known methods for identifying Du Bois singularities, and provide some new ones.  Steenbrink's following result is one of the most useful.
\begin{thm}[{\cite[Thm 3]{St81}}]\label{T:dbnorm}  Let $X$ be a variety, with $\pi:\widetilde{X}\to X$ its normalization and $\sC=\ann_{\sO_X}(\pi_{*}\sO_{\widetilde{X}}/\sO_{X})$ the conductor ideal sheaf of the map $\pi$.  Define $\Sigma\sub X$ to be the subvariety defined by $\sC$, and let $E=\pi^{-1}(\Sigma)$.  Suppose $\widetilde{X}, E$ and $\Sigma$ all have Du Bois singularities.  Then $X$ has Du Bois singularities.
\end{thm}

\begin{ex}\label{E:db1}
The pinch point is a Du Bois singularity.  The
normalization of its coordinate ring is $\C[y_1,y_2,x_4,\ldots,x_n]$, which defines a smooth (and therefore Du Bois) scheme.  The normalization map is given by
\begin{align*}
x_1 &\mapsto y_1 y_2 \\
x_2 &\mapsto y_2 \\
x_3 &\mapsto y_1^2 \\
x_i &\mapsto x_i, \, i\geq 4.
\end{align*}
The conductor is the ideal $(x_1,x_2)$, which defines a smooth subscheme.
To obtain the preimage of this subscheme, we take the image of the
conductor in the normalization, which is the ideal $(y_2)$; thus we can apply the theorem.
\end{ex}

\begin{ex}\label{E:db2}
A double normal crossing singularity is also Du Bois.  The normalization of its coordinate ring is $\C[x_1,\ldots,x_n]/(x_1)\oplus \C[x_1,\ldots,x_n]/(x_2)$, which defines a smooth variety.  The conductor is given by the ideal $(x_1, x_2)$, which defines a smooth subscheme; the preimage of this subscheme is the direct sum of two copies of $\C[x_1,\ldots,x_n]/(x_1, x_2)$, which is again smooth.
\end{ex}

Pinch points and double normal crossing points are the only singularities of semismooth schemes, so the above examples give the following.
\begin{prop} If $X$ is semismooth, then $X$ is Du Bois.
\end{prop}

Karl Schwede recently proved the following interesting connection between F-injective and Du Bois singularities (see \cite{Fed83} and many others for details on F-injective and F-pure singularities).
\begin{thm}[{\cite[Thm. 6.4.3]{KS06}}]\label{T:fpdb} Suppose that $X$ is a variety over $\C$ of dense F-injective type.  Then $X$ has Du Bois singularities.
\end{thm}
Combined with \cite[Prop 2.1]{Fed83} and \cite[Lemma 3.3]{Fed83}, this gives the corollary we shall use:

\begin{cor}\label{C:fpdb} Suppose $X=V(f)$ is a variety defined (locally) by a single equation $f\in \C[x_1,\ldots,x_n]$.  Assume $P=(0,\ldots,0)\in X$ is a singular point.  If $f^{p-1}\notin (x_1^p,\ldots,x_n^p)$ for all but finitely many primes $p$, then $P\in X$ is a Du Bois singularity.
\end{cor}

\subsection{Products of Du Bois Schemes}

Does a product of schemes with Du Bois singularities also have Du Bois singularities?  This is a natural question to ask, but seems to be missing from the scant literature on this class of singularities.  Such a ``product theorem'' is the final piece required before proceeding to the proof of the Main Theorem.

One more type of singularity plays a role in the product theorem.

\begin{defn} If $X$ is a variety, we say that $X$ has \emph{generalized simple normal crossings} if for each singular point $x\in X$ we have an analytic isomorphism
\[
\widehat{\sO}_{X,x}\cong \C[[x_1,\ldots,x_n]]/(I_1\cap\cdots\cap I_k),
\]
where each $I_k$ is generated by coordinate functions, i.e., $I_k=(x_{k_1},\ldots,x_{k_j})$.
\end{defn}

\begin{lemma}\label{L:gsnc} Suppose that $X$ has generalized simple normal crossings.  Then $X$ has Du Bois singularities.
\end{lemma}

\begin{proof}
Note the following ``gluing'' fact of Du Bois singularities: if $X$ is a variety over $\C$ with components $X_1$ and $X_2$ such that $X_1, X_2$ and $X_1\cap X_2$ have Du Bois singularities, then $X$ has Du Bois singularities (\cite[Thm. 5.2.1]{KS06}).

We proceed by induction on the maximum dimension of a component and the number of components.  If $X$ is the union of two components of any dimension, it has Du Bois singularities by the gluing property (here the intersection is actually smooth).  Now assume that the conclusion holds when all components are less than dimension $d$, and that it holds for $n$ components of dimension $d$.  Suppose that we have $n+1$ components of dimension at most $d$.  Then $X$ can be expressed as the union of the components $X_1$ and $X_2$, where $X_1$ is the union of the first $n$ components and $X_2$ is the last component.  By the induction hypothesis $X_1$ is Du Bois, and $X_2$ is smooth.  $X_1\cap X_2$ has generalized simple normal crossing singularities and components of dimension $d-1$ or less, so by the induction hypothesis it also has Du Bois singularities.  By the gluing property, we conclude that $X$ has Du Bois singularities.
\end{proof}

\begin{lemma}\label{L:prod} Suppose $X$ and $Y$ have generalized simple normal crossings.  Then $X\times Y$ also has generalized simple normal crossings.
\end{lemma}

\begin{proof}  Let $z=(x,y)\in X\times Y$.  By assumption, $\widehat{\sO}_{X,x}\cong \C[[x_1,\ldots,x_n]]/I$ and $\widehat{\sO}_{Y,y}\cong \C[[y_1,\ldots,y_m]]/J$ where $I$ and $J$ are intersections of ideals generated by coordinate functions.  Thus
\begin{align*}
\widehat{\sO}_{X\times Y,z}&\cong \C[[x_1,\ldots,x_n]]/I\tensor_{k}\C[[y_1,\ldots,y_m]]/J\\
&\cong \C[[x_1,\ldots,x_n,y_1,\ldots,y_m]]/IJ\\
&=\C[[x_1,\ldots,x_n,y_1,\ldots,y_m]]/(I\cap J),
\end{align*}
where the last equality follows from the fact that $I$ and $J$ are ideals in disjoint polynomial rings.  Thus $X\times Y$ has generalized simple normal crossings.
\end{proof}

\begin{thm}\label{T:Prod} Suppose $X_1$ and $X_2$ are varieties (over $\C$) which have Du Bois singularities.  Then $X_1\times X_2$ also has Du Bois singularities.
\end{thm}
\begin{proof}
The property of having Du Bois singularities is local (in the \'etale topology), so by restricting to affine open sets we may assume that $X_1$ and $X_2$ are affine.  Embed each $X_i$ into a smooth variety $Y_i$, and let $f_i:\widetilde{Y}_i\to Y_i$ be a strong log resolution of $X_i$.  Denote by $E_i$ the (reduced) pre-image of $X_i$.  By \cite[Thm. 5.3.4]{KS06}, we have $\DB{0}{X_i}\qis Rf_{i*}\sO_{E_i}$.  Since $X_i$ was assumed to have Du Bois singularities, we have $\sO_{X_i}\qis Rf_{i*}\sO_{E_i}$.

Note that $\widetilde{Y}_1\times\widetilde{Y}_2$ is smooth (hence has rational singularities) and $g$ is an isomorphism outside of $X_1\times X_2$.  By assumption, $E_1$ and $E_2$ have generalized simple normal crossings.  Hence $E_1\times E_2$ has generalized simple normal crossings by Lemma \ref{L:prod}, and so has Du Bois singularities by Lemma \ref{L:gsnc}.  In particular, this means that $Rg_*\sO_{E_1\times E_2}\qis\DB{0}{X_1\times X_2}$ by \cite[Thm. 5.3.4]{KS06}.

Now consider the following diagram:
\begin{equation*}\begin{CD}
E_1 @<{\widetilde{\pi}_1}<< E_1\times E_2 @. @.\\
@V{f_1}VV @VV{g_1}V @. @. \\
X_1 @<{\pi_1}<< X_1\times E_2 @>{\widetilde{\pi}_2}>> E_2 \\
@. @VV{g_2}V @VV{f_2}V\\
@. X_1\times X_2 @>>{\pi_2}> X_2
\end{CD}
\end{equation*}
(Note that $g_1$ and $g_2$ are the obvious maps such that $g_2\circ g_1=g|_{E_1\times E_2}$, and $\widetilde{\pi}_i$ and $\pi_i$ are the usual projection maps.)  Applying \cite[Thm. III.9.3]{Ha77} to the upper left square in the diagram, we find that $\pi_1^*Rf_{1*}\sO_{E_1}\cong Rg_{1*}(\widetilde{\pi}_1^*\sO_{E_1})$.  Since $\widetilde{\pi}_1$ is projection onto a factor, we have $\widetilde{\pi}_1^*\sO_{E_1}\cong\sO_{E_1\times E_2}$; by assumption, $Rf_{1*}\sO_{E_1}\qis\sO_{X_1}$, so we obtain $Rg_{1*}(\sO_{E_1\times E_2})\qis\pi_1^*\sO_{X_1}\cong\sO_{X_1\times E_2}$ (the last isomorphism is again due to the fact that $\pi_1$ is a projection map).  Applying this same argument again to the lower right square gives $Rg_{2*}\sO_{X_1\times E_2}\qis\sO_{X_1\times X_2}$. Thus we see that $Rg_{*}\sO_{E_1\times E_2}\qis\sO_{X_1\times X_2}$; combined with the above, we see that $\sO_{X_1\times X_2}\qis\DB{0}{X_1\times X_2}$, so that $X_1\times X_2$ is Du Bois.
\end{proof}

\section{Singularities of Generic Projections}

\subsection{Proof of the Main Theorem}

Recall that the set of linear projections $\P^m\to\P^r$ is in bijection with the $(m-r)$-dimensional linear subspaces of $\P^m$, and can thus be identified with the closed points of the Grassmannian variety $G(m, m-r-1)$.  We say that a \emph{generic projection} has property $P$ if the collection of points with property $P$ forms an open dense subset of the Grassmannian.  With this in mind, we proceed to the proof of the main theorem.

\begin{mthm} Let $Y\subset\P^N$ be a smooth projective variety of dimension $n$, $n\leq5$, embedded via the $d$-uple embedding with $d\geq 3n$.  Let $X\subset\P^{n+1}$ be the image of $Y$ under a generic projection $\pi:Y\to\P^{n+1}$.  Then $X$ has Du Bois singularities.
\end{mthm}

\begin{proof}
\cite[\S 13.2]{Rob75} provides a list of the possible local analytic isomorphism classes of the singularities that arise from such a generic projection.  Since the Du Bois complex is local in the \'etale toplogy, it will suffice to show that in each case these define Du Bois singularities.

\begin{enumerate}[\upshape {Case} (1)]
\item[Case (0)] Let $R=\C[[x_1,\ldots,x_{n+1}]]/(x_1\cdots x_d)$. In this case, we have $f_1=x_1\cdots x_d$.  Reducing to characteristic $p$, we see that $f_1^{p-1}=x_1^{p-1}\cdots x_n^{p-1}$, and $\mathfrak{m}^{[p]}=(x_1^p,\ldots,x_{n+1}^p)$.  We clearly have $f_1^{p-1}\notin\mathfrak{m}^{[p]}$ for all primes $p$, so $R$ has Du Bois singularities by \ref{C:fpdb}.

\item[Case (1a)] The pinch point was shown to be Du Bois in Example \ref{E:db1}.

\item[Case (1b)] Let $R=\C[[x_1,\ldots,x_{n+1}]]/(x_{n}^3+\Phi_4+\Phi_5)$, with
\begin{align*}
\Phi_4 &= x_{1}^2 x_{3} x_{n}-x_{1}^3x_{n+1}+2x_{2}x_{3}x_{n}^2-3x_{1}x_{2}x_{n}x_{n+1},\\
\Phi_5 &= x_{2}^2x_{3}^2x_{n}-x_{1}x_{2}^2x_{3}x_{n+1}-x_{2}^3x_{n+1}^2;
\end{align*}
denote the full polynomial generating the ideal by $f_3$.  Examination of the monomials occurring in $f_3$ shows that there is a term of the form $(-3x_1x_2x_nx_{n+1})^{p-1}$ in $f_3^{p-1}$.  No other product of monomials in $f_3$ can generate a monomial of the form $(x_1x_2x_nx_{n+1})^k$, so the coefficient of $(x_1x_2x_nx_{n+1})^{p-1}$ is nonzero for $p\neq3$.  Since this monomial is not in $\mathfrak{m}^{[p]}$, it follows that $f_3\notin\mathfrak{m}^{[p]}$.  Again, \ref{C:fpdb} implies this point id Du Bois.

\item[Case (2a)] The ring $R=\C[[x_1,\ldots,x_{n+1}]]/(x_1(x_{n}^2-x_{2}^2x_{n+1}))$ is the coordinate ring of a product $X_1\times X_2$ (where we take $X_1$ to be a pinch point and $X_2$ to be a hyperplane), hence defines a Du Bois singularity by Theorem \ref{T:Prod}.

\item[Case (2b)] $R=\C[[x_1,\ldots,x_{n+1}]]/(x_1(x_{n}^3+\Psi_4+\Psi_5))$,
where
\[
\Psi_i = \Phi_i(x_2,x_3,x_4,x_n,x_{n+1}),
\]
is also Du Bois by Theorem \ref{T:Prod} and case (1b).

\item[Case (2c)] $R=\C[[x_1,\ldots,x_{n+1}]]/((x_n^2-x_1^2x_{n+1})(x_{n-2}^2-x_2^2x_{n-1}))$ defines a Du Bois singularity by the product theorem and case (1a).

\item[Case (3)] $R=\C[[x_1,\ldots,x_{n+1}]]/(x_1x_2(x_{n}^2-x_{3}^2x_{n+1}))$ is Du Bois by the product theorem and case (2a).

\item[Case (4)] Finally, $R=\C[[x_1,\ldots,x_{n+1}]]/(x_1x_2x_3(x_{n}^2-x_{4}^2x_{n+1}))$ is Du Bois by the product theorem and case (3).
\end{enumerate}
\end{proof}

\subsection{Du Bois versus Semi Log Canonical}

With a few additional results, we can use the Main Theorem to address semi log canonicity of generic projection hypersurfaces.  Toward that end, we demonstrate that Du Bois singularities are in fact semi log canonical under certain additional assumptions.

\begin{lemma}\label{L:dbres}
Let $f:Y\to X$ be a good semiresolution. Assume that $\pi$ is an isomorphism
outside a closed subscheme $\Sigma\subset X$, with preimage $f^{-1}(\Sigma)=E$.  Then
there exists an exact triangle
\[
\DB{0}{X}\to\DB{0}{\Sigma}\oplus Rf_{*}\sO_{Y}\to Rf_{*}\sO_{E}\xrightarrow{\text{\upshape +1}}.
\]
\end{lemma}

\begin{proof}
Since $Y$ is semismooth, $\DB{0}{Y}\qis\sO_{Y}$.  The assumption that $f$ is a good semiresolution means that $E$ has only normal crossing singularities, so that $\DB{0}{E}\qis\sO_{E}$.  The result then follows from \cite[Prop. 4.11]{DB81}.
\end{proof}

The following theorem generalizes a similar result for log canonical singularities, due to S\'andor Kov\'acs (\cite[Thm 3.6]{SK99}).

\begin{thm}\label{T:slcdb} Let $X$ be an $S_2$ scheme which is semismooth in codimension one, and assume that $K_{X}$ is Cartier and $X$ has Du Bois singularities.  Then $X$ is semi log canonical.
\end{thm}

\begin{proof} Let $f:Y\to X$ be a good semiresolution of $X$, with $W\sub X$ the set outside which $f$ is an isomorphism, and $E=f^{-1}(S)$.  There exists a natural morphism $\phi:Rf_{*}\sO_{Y}(-E)\to Rf_{*}\sO_{Y}$ arising from the short exact sequence
\[
0\to\sO_{Y}(-E)\to\sO_{Y}\to\sO_{E}\to 0.
\]
We note that $Rf_{*}\sO_{Y}(-E)\to Rf_{*}\sO_{E}$ is the zero map, from which it follows (via the exact triangle in Proposition \ref{L:dbres}) that $\phi$ factors through $\DB{0}{X}$.  By assumption, $\DB{0}{X}\qis\sO_{X}$, so we obtain a morphism $Rf_{*}\sO_{Y}(-E)\to\sO_{X}$ which is a quasi-isomorphism on $X\backslash S$.  Applying $R\sHom_{X}(-,\omega^{\cdot}_{X})$, we obtain a morphism $\omega^{\cdot}_X\to R\sHom_{X}(Rf_{*}\sO_{Y}(-E),\omega^{\cdot}_{X})$.  But then we have
\begin{align*}
R\sHom_{X}(Rf_{*}\sO_{Y}(-E),\omega^{\cdot}_{X})&\cong Rf_{*}R\sHom_{Y}(\sO_{Y}(-E),\omega^{\cdot}_{Y})\\
&\cong Rf_{*}\left(R\sHom_{Y}(\sO_{Y}(-E),\sO_{Y})\tensor\omega^{\cdot}_{Y}\right),
\end{align*}
where the first isomorphism follows from Grothendieck duality, and the second follows from the fact that $Y$ is Gorenstein.  Since $E$ is a Cohen-Macaulay divisor (it has only simple normal crossings), this last term is isomorphic to $Rf_{*}\omega_{Y}(E)[n]$, so in fact we have a morphism $\omega^{\cdot}_{X}\to Rf_{*}\omega_{Y}(E)[n]$.  Taking the $-n$th cohomology gives a morphism $\omega_{X}\to f_{*}\omega_{Y}(E)$ which is an isomorphism on $X\backslash S$.  Adjointness produces a nonzero morphism $f^{*}\omega_{X}\to\omega_{Y}(E)$.  Since
$f^{*}\omega_{X}$ is a line bundle this implies $f^{*}\omega_{X}\subseteq\omega_{Y}(E)$, whence $X$ is semi log canonical.
\end{proof}

\begin{cor}\label{C:DB} Let $X$ be a seminormal Gorenstein scheme with Du Bois singularities.  Then $X$ is semi log canonical.
\end{cor}

\begin{proof} The Gorenstein assumption implies that $K_X$ is Cartier, and also that $X$ is G1, i.e., $X$ is $S_2$ and Gorenstein in codimension 1.  By Theorem \cite[Thm 9.10]{GT80}, a seminormal G1 scheme is semismooth in codimension one. The result then follows immediately from the theorem.
\end{proof}

\begin{cor} Let $X$ be as in the Main Theorem.  Then $X$ has semi log canonical singularities.
\end{cor}

\begin{proof} Since $X$ is a hypersurface, it is a complete intersection; in particular,
$X$ is Gorenstein.  By \cite[Thm 3.7]{GT80}, $X$ is seminormal, so \ref{C:DB} implies
that $X$ is semi log canonical.
\end{proof}

\subsection{Counterexamples in Higher Dimensions}

Ideally we would like a more general result that generic projection hypersurfaces are \emph{always} Du Bois (and thus semi log canonical).  The following results illustrate that counterexamples exist, though all known examples have high dimension.

\begin{thm}\label{T:notslc} Suppose $X\subset\P^{n+1}$ is a hypersurface, and that there exists some point $x\in X$ having multiplicity $\mu>n+1$.  Then $X$ is not semi log canonical.
\end{thm}

\begin{proof}
Let $f:Z\to\P^{n+1}$ be the blow-up at $x$, and let $g:X'\to X$ be the restriction of $f$ to the strict transform of $X$.  Since $f$ is a blow-up at a point, we have (by \cite[Ex. II.8.5]{Ha77})
\[
K_Z\equiv f^{*}K_{\P^{n+1}} + nE.
\]

Similarly, by the definition of $X'$ and our choice of $x$ we have
\[
X'=g^{*}X-\mu E|_{X'},
\]
where we are abusing notation and identifying $X$ and $X'$ with the
corresponding divisors.  Applying the adjunction formula, we obtain
\begin{align*}
K_{X'}&\equiv (K_Z+X')|_{X'}\\
&\equiv (f^{*}K_{\P} + X' + nE)|_{X'}\\
&\equiv (f^{*}K_{\P}+g^{*}X-\mu E|_{X'} + nE)|_{X'}\\
&\equiv f^{*}(K_{\P}+X)+(n-\mu)E|_{X'}\\
&\equiv g^{*}K_{X}+(n-\mu)E|_{X'}.
\end{align*}
Note that a good semiresolution of $X'$ also produces a good semiresolution of $X$.  Furthermore, semiresolving $X'$ will not increase the coefficient of $E$.  Since $n-\mu<-1$, we conclude that $X$ is not semi log canonical.
\end{proof}

Given a morphism $f:Y\to X$, we denote by $S_i(f)$ the locus where the induced morphism on tangent spaces $df$ drops rank by $i$, i.e.,
\[
S_i(f)=\{y\in Y\mid\rank df_{y}\leq\dim Y-i\}.
\]

\begin{prop}\label{P:mult} Let $f:Y\to X\subset\P^{n+1}$ be a finite morphism, and suppose $y\in S_i(f)$.  Then the point $f(y)\in X$ has multiplicity at least $2^i$.
\end{prop}

\begin{proof} Since $y\in S_i(f)$, the map $df_y:T_yY\to T_{f(y)}X$ has rank at most $n-i$.  To compute the multiplicity of $f(y)$ on $X$, we compute the intersection multiplicity of a general line $L$ with $X$ at $f(y)$.  $L$ is determined by $n$ linear forms, say $l_1,\ldots, l_n$.  We can compute the intersection multiplicity by pulling back the $l_i$ to $Y$, where they generate $n$ hypersurfaces.  Since $df$ drops rank by $i$ at $y$, no more than $n-i$ of the equations defining these hypersurfaces have independent linear terms at $y$.  Without loss of generality, we may assume that the remaining $i$ equations have at least degree 2 at $y$; thus $(f^{*}L).Y$ has multiplicity at least $2^i$ at $y$.  Therefore $f(x)$ also has multiplicity at least $2^i$.
\end{proof}

\begin{cor} Let $X\subset\P^{31}$ be a generic projection hypersurface obtained via $\pi:Y\to\P^{31}$, where $\OM^1_Y$ is nef.  Then $X$ is not semi log canonical.
\end{cor}

\begin{proof} The hypothesis that $\OM^1_Y$ is nef, together with 
\cite[Cor. 7.2.18]{Laz1}, implies that $S_5(\pi)\neq\emptyset$.  Proposition \ref{P:mult} implies that for any $y\in S_5(\pi)$, the image $f(y)$ has multiplicity at least $2^5=32$.  The result then follows by Theorem \ref{T:notslc}.
\end{proof}

\begin{ex} If $Y$ is a smooth projective scheme over $\C$ which is uniformized by $\B^n\subset\C^n$, then $\OM^1_Y$ is ample, and thus nef (cf. \cite[6.3.36]{Laz2}, \cite{ZL86}).
\end{ex}

\begin{ex} Let $Y_1,\ldots,Y_m$ be smooth projective varieties of dimension $d\geq1$, each with big cotangent bundle.  If $Y\subset Y_1\times\cdots\times Y_m$ is a general linear section, with
\[
\dim Y\leq\frac{d(m+1)+1}{2(d+1)},
\]
then $\OM^1_Y$ is ample (\cite{D05}).
\end{ex}

\begin{ex} If $Y$ is the complete intersection of at least $n/2$ sufficiently ample general hypersurfaces in an abelian variety of dimension $n$, then $\OM^1_Y$ is ample (\cite{D05}).
\end{ex}

\begin{ex} If $Y$ is a projective variety over $\C$ whose universal covering space is a bounded domain in $\C^n$ or a Stein manifold, then $\OM^1_Y$ is nef (\cite{K97}).
\end{ex}

\bibliographystyle{skalpha}
\bibliography{biblio}

\end{document}

%% file: hart.tex
\DeclareFontFamily{OMS}{rsfs}{\skewchar\font'60}
\DeclareFontShape{OMS}{rsfs}{m}{n}{<-5>rsfs5 <5-7>rsfs7 <7->rsfs10 }{}
\DeclareSymbolFont{rsfs}{OMS}{rsfs}{m}{n}
\DeclareSymbolFontAlphabet{\scr}{rsfs}

\newcommand{\sC}{\scr{C}}

\newcommand{\sO}{\scr{O}}

\newcommand{\OM}{\Omega}

\newcommand{\qis}{\simeq_{\text{qis}}}

\newcommand{\ul}{\underline}

\newcommand{\B}{\mathfrak B}

\newcommand{\C}{\mathbb {C}}

\newcommand{\Q}{\mathbb {Q}}

\renewcommand{\P}{\mathbb{P}}

\newcommand{\sub}{\subseteq}
\newcommand{\tensor}{\otimes}

\DeclareMathOperator{\ann}{{Ann}}

\DeclareMathOperator{\rank}{{rank}}

\DeclareMathOperator{\Gr}{{Gr}}

\DeclareMathOperator{\sHom}{\scr{H}\textit{om}}

\def\spec#1.#2.{{\bold S\bold p\bold e\bold c}_{#1}#2}
\def\proj#1.#2.{{\bold P\bold r\bold o\bold j}_{#1}\sum #2}
\def\ring#1.{\scr O_{#1}}
\def\map#1.#2.{#1 \to #2}
\def\longmap#1.#2.{#1 \longrightarrow #2}
\def\factor#1.#2.{\left. \raise 2pt\hbox{$#1$} \right/
\hskip -2pt\raise -2pt\hbox{$#2$}}
\def\pe#1.{\mathbb P(#1)}
\def\pr#1.{\mathbb P^{#1}}

\newcommand{\DB}[2]{{\ul\Omega}^{#1}_{#2}}

\def\coh#1.#2.#3.{H^{#1}(#2,#3)}
\def\ccoh#1.#2.#3.{\check{H}^{#1}(#2,#3)}
\def\dimcoh#1.#2.#3.{h^{#1}(#2,#3)}
\def\hypcoh#1.#2.#3.{\mathbb H_{\vphantom{l}}^{#1}(#2,#3)}
\def\loccoh#1.#2.#3.#4.{H^{#1}_{#2}(#3,#4)}
\def\dimloccoh#1.#2.#3.#4.{h^{#1}_{#2}(#3,#4)}
\def\lochypcoh#1.#2.#3.#4.{\mathbb H^{#1}_{#2}(#3,#4)}
\def\ses#1.#2.#3.{0  \longrightarrow  #1   \longrightarrow 
 #2 \longrightarrow #3 \longrightarrow 0} 
\def\sesshort#1.#2.#3.{0
 \rightarrow #1 \rightarrow #2 \rightarrow #3 \rightarrow 0}
\def\iff#1#2#3{
    \hfil\hbox{\hsize =#1 \vtop{\noin #2} \hskip.5cm
    \lower.5\baselineskip\hbox{$\Leftrightarrow$}\hskip.5cm
    \vtop{\noin #3}}\hfil\medskip}
\def\myoplus#1.#2.{\underset #1 \to {\overset #2 \to \oplus}}

%% file: slcgenproj.bbl
\newcommand{\SortNoop}[1]{}
\providecommand{\bysame}{\leavevmode\hbox to3em{\hrulefill}\thinspace}
\providecommand{\MR}{\relax\ifhmode\unskip\space\fi MR}
\providecommand{\MRhref}[2]{%
  \href{http://www.ams.org/mathscinet-getitem?mr=#1}{#2}
}
\providecommand{\href}[2]{#2}
\begin{thebibliography}{Laz04b}

\bibitem[Deb05]{D05}
{\sc O.~Debarre}: \emph{Varieties with ample cotangent bundle}, Compos. Math.
  \textbf{141} (2005), no.~6, 1445--1459. {\sf\scriptsize MR2188444}

\bibitem[DB81]{DB81}
{\sc P.~Du~Bois}: \emph{Complexe de de {R}ham filtr\'e d'une vari\'et\'e
  singuli\`ere}, Bull. Soc. Math. France \textbf{109} (1981), no.~1, 41--81.
  {\sf\scriptsize MR613848 (82j:14006)}

\bibitem[Fed83]{Fed83}
{\sc R.~Fedder}: \emph{{$F$}-purity and rational singularity}, Trans. Amer.
  Math. Soc. \textbf{278} (1983), no.~2, 461--480. {\sf\scriptsize MR701505
  (84h:13031)}

\bibitem[GT80]{GT80}
{\sc S.~Greco and C.~Traverso}: \emph{On seminormal schemes}, Compositio Math.
  \textbf{40} (1980), no.~3, 325--365. {\sf\scriptsize 81j:14030}

\bibitem[Har77]{Ha77}
{\sc R.~Hartshorne}: \emph{Algebraic geometry}, Springer-Verlag, New York,
  1977. {\sf\scriptsize 57 \#3116}

\bibitem[Kol90]{K90}
{\sc J.~Koll{\'a}r}: \emph{Projectivity of complete moduli}, J. Differential
  Geom. \textbf{32} (1990), no.~1, 235--268. {\sf\scriptsize MR1064874
  (92e:14008)}

\bibitem[Kov99]{SK99}
{\sc S.~Kov{\'a}cs}: \emph{Rational, log canonical, {D}u {B}ois singularities:
  on the conjectures of {K}oll\'ar and {S}teenbrink}, Compositio Math.
  \textbf{118} (1999), no.~2, 123--133. {\sf\scriptsize MR1713307
  (2001g:14022)}

\bibitem[Kra97]{K97}
{\sc H.~Kratz}: \emph{Compact complex manifolds with numerically effective
  cotangent bundles}, Doc. Math. \textbf{2} (1997), 183--193 (electronic).
  {\sf\scriptsize MR1464070 (98j:32033)}

\bibitem[Laz04a]{Laz1}
{\sc R.~Lazarsfeld}: \emph{Positivity in algebraic geometry. {I}}, Ergebnisse
  der Mathematik und ihrer Grenzgebiete. 3. Folge. A Series of Modern Surveys
  in Mathematics [Results in Mathematics and Related Areas. 3rd Series. A
  Series of Modern Surveys in Mathematics], vol.~48, Springer-Verlag, Berlin,
  2004, Classical setting: line bundles and linear series. {\sf\scriptsize
  MR2095471 (2005k:14001a)}

\bibitem[Laz04b]{Laz2}
{\sc R.~Lazarsfeld}: \emph{Positivity in algebraic geometry. {II}}, Ergebnisse
  der Mathematik und ihrer Grenzgebiete. 3. Folge. A Series of Modern Surveys
  in Mathematics [Results in Mathematics and Related Areas. 3rd Series. A
  Series of Modern Surveys in Mathematics], vol.~49, Springer-Verlag, Berlin,
  2004, Positivity for vector bundles, and multiplier ideals. {\sf\scriptsize
  MR2095472 (2005k:14001b)}

\bibitem[Rob71]{Rob71}
{\sc J.~Roberts}: \emph{Generic projections of algebraic varieties}, Amer. J.
  Math. \textbf{93} (1971), 191--214. {\sf\scriptsize MR0277530 (43 \#3263)}

\bibitem[Rob75]{Rob75}
{\sc J.~Roberts}: \emph{Singularity subschemes and generic projections}, Trans.
  Amer. Math. Soc. \textbf{212} (1975), 229--268. {\sf\scriptsize MR0422274 (54
  \#10265)}

\bibitem[Sch06]{KS06}
{\sc K.~Schwede}: \emph{On f-injective and du bois singularities}, Ph.D.
  thesis, University of Washington, 2006.

\bibitem[Ste81]{St81}
{\sc J.~H.~M. Steenbrink}: \emph{Cohomologically insignificant degenerations},
  Compositio Math. \textbf{42} (1980/81), no.~3, 315--320. {\sf\scriptsize
  MR607373 (84g:14011)}

\bibitem[ZL86]{ZL86}
{\sc M.~G. Za{\u\i}denberg and V.~Y. Lin}: \emph{Finiteness theorems for
  holomorphic mappings}, Current problems in mathematics. Fundamental
  directions, Vol.\ 9 (Russian), Itogi Nauki i Tekhniki, Akad. Nauk SSSR
  Vsesoyuz. Inst. Nauchn. i Tekhn. Inform., Moscow, 1986, pp.~127--193, 272.
  {\sf\scriptsize MR860611}

\end{thebibliography}
